%% file: main.tex
 \documentclass[tablecaption=bottom,wcp]{jmlr} 

\usepackage{booktabs}
\usepackage[load-configurations=version-1]{siunitx} 
\usepackage{hyperref}

\usepackage{subcaption}

\theorembodyfont{\upshape}
\theoremheaderfont{\scshape}
\theorempostheader{:}
\theoremsep{\newline}

\jmlrvolume{1}
\jmlryear{2018}
\jmlrsubmitted{19 March 2018}
\jmlrpublished{publication date}
\jmlrworkshop{Submitted to COPA 2018} 

\title[Misspecification GP interpolation error]{Interpolation error of misspecified Gaussian process regression}

  \author{\Name{A. Zaytsev} \Email{a.zaytsev@skoltech.ru} \\
   \Name{E. Romanenkova} \Email{shulgina@phystech.edu}\\
   \Name{D. Ermilov} \Email{dmitrii.ermilov@skoltech.ru}\\
   \addr Skoltech, IITP RAS, MIPT}

\input{mydef.tex}

\renewcommand{\bbR}{\mathbb{R}}
\newcommand{\wtheta}{\theta'}

\editor{Editor's name}

\begin{document}

\maketitle

\begin{abstract}

An interpolation error is an integral of the squared error of a regression model over a domain of interest.
We consider the interpolation error for the case of misspecified Gaussian process regression: used covariance function differs from the true one.
We derive the interpolation error for an infinite grid design of experiments.
In particular, we show that for $\mathrm{Matern}_{\frac{1}{2}}$ covariance function poor estimation of parameters only slightly affects the quality of interpolation. 
Then we proceed to numerical experiments that consider the misspecification
for the most common covariance functions including other Matern and squared exponential covariance functions.
For them, the quality of estimates of parameters affects the interpolation error.
\end{abstract}

\section{Introduction}
\label{sec:intro}

Gaussian process regression or kriging is widely used 
for construction of regression models~\cite{rasmussen2006,burnaev2016regression,cressie2015statistics}.
The main assumption of these approaches is that the target function 
is a realization of Gaussian process model with a given spectral density (or equivalently a covariance function) and mean functions.

For each approach it is of great importance to get a measure of quality of a regression model.
Popular choice in literature is an interpolation error~\cite{golubev13interpolation,le2015asymptotic}: 
an expected squared error of interpolation integrated over a domain of interest for a given approach of a regression model construction.

There a number of problem statements relevant to this general problem.
Classical approaches imply that the true model is known and coincides with the one used for construction of a regression model~\cite{stein2012interpolation}.
Modern approaches more often consider a minimax problem statement ~\cite{zaytsev2017minimax} 
or a misspecified problem statement~\cite{vaart2011information}. 
In the minimax problem statement we assume that the true model belongs to a certain class of models and try to find the interpolation error in the worst case~\cite{golubev13interpolation}.
In the misspecified problem statement we specify how a used models differs from the true models~\cite{panov2016nonasymptotic}.

Let us elaborate in more details the misspecified problem statement for Gaussian process regression.
For real problems one doesn't know the true Gaussian process regression model, while The usual assumption is that the spectral density (a Fourier transform of the covariance function) belongs to a given parametric family and the mean value is zero.
After selection of a parametric family one estimates parameters of a spectral density using approaches similar to the maximum likelihood approach or Bayesian approach~\cite{zaytsev2014properties}
Quality of estimation of parameters varies~\cite{zaytsev2014properties,bachoc2018asymptotic}.
Moreover, smoothness of the target function is often unknown.
So it is hard to select a parametric family of spectral densities.
Thus, bad estimates of spectral density and wrong choice of a parametric family lead to a difference between the true regression 
model and the used regression model.

Our goal is to obtain exact expression for interpolation error in a misspecification case. 
The assumptions are similar to used in the state of the art:
Gaussian process is stationary,
the design of experiments is an infinite grid with a given step along each dimension.
The grid designs of experiments are often used due to their low computational complexities~\cite{belyaev2015gaussian}.
Moreover, numerical experiments show that these assumption don't significantly affect the results~\cite{zaytsev2017minimax}.
Using obtained expression as a tool we are able to consider 
widely used setups for Gaussian process regression taking into account possible model misspecification.
We consider the squared exponential function and the Matern covariance functions with $\nu = \frac12$~\cite{minasny2005matern}.

The article has the following sections:
\begin{itemize}
\item Section~\ref{sec:prior} describes the prior results in this area in more details;
\item Section~\ref{sec:known} describes results for usage of known covariance function and minimax case;
\item Section~\ref{sec:wrong} describes results for the case when the true covariance function differs from the used one and examines in more details the case of model misspecification for the Matern covariance function;
\item Section~\ref{sec:experiments} contains results of numerical experiments;
\end{itemize}

\section{Related work}
\label{sec:prior}

Classical approaches imply that the true model is known and coincides with the one used for construction of a regression model.
The first results in this area go back to \cite{kolmogorov1992interpolation} and \cite{wiener1949extrapolation}.
A.Kolmogorov and N.Wiener simultaneously obtained mean squared errors at a point in an interpolation and an extrapolation problem statements
with all training points lying on a grid.
An article~\cite{le2015asymptotic} considered the integrated mean squared interpolation error for a Gaussian process with noise if the sample size tends to infinity.

Modern approaches more often consider a minimax problem statement.
An article~\cite{golubev13interpolation} considered the minimax interpolation error for a Sobolev class of Gaussian processes for a segment if the training sample is an infinite grid.
More recent article~\cite{zaytsev2017minimax} considered multivariate scenario, while considering Gaussian processes with an upper bound only for a sum of squares of the first partial derivatives of Gaussian process realization.

Another branch of modern results considers a misspecified problem statement.
For a review of results for a squared error at a single point see book~\cite{stein2012interpolation}
More general papers~\cite{,van2008rates,vaart2011information} consider the case of mean squared error for an area, while their results are not directly applicable in a practice-related problems due to complex assumptions.
Note also that these articles as well as~\cite{suzuki2012pac,castillo2008lower} provide upper bound.
An article \cite{bachoc2013cross} considered empirical comparison of the interpolation error for cross validation and maximum likelihood estimates, 
while theoretical properties of these approaches are investigated in
more details in~\cite{bachoc2018asymptotic}, 
while the focus is not on the interpolation error itself, but on the quality of parameter estimation.

\section{Interpolation error and minimax interpolation error}
\label{sec:known}

Let us introduce interpolation for the case with no misspecification  and the minimax case.
All results in this section are provided in a way similar to~\cite{zaytsev2017minimax}. 

For $\bbR^{\iD}$ there is a stationary Gaussian process $f(\cdot)$ with the covariance function
$R(\vecX)$.
The spectral density~\cite{stein2012interpolation}  is defined as 
\[
F(\vecO) = \frac{1}{2 \pi} \int_{\bbR^{\iD}} e^{-\mathrm{i} \vecO \vecX} R(\vecX) d\vecX.
\]

\begin{figure}
   \centering
   \includegraphics[width=0.4\textwidth]{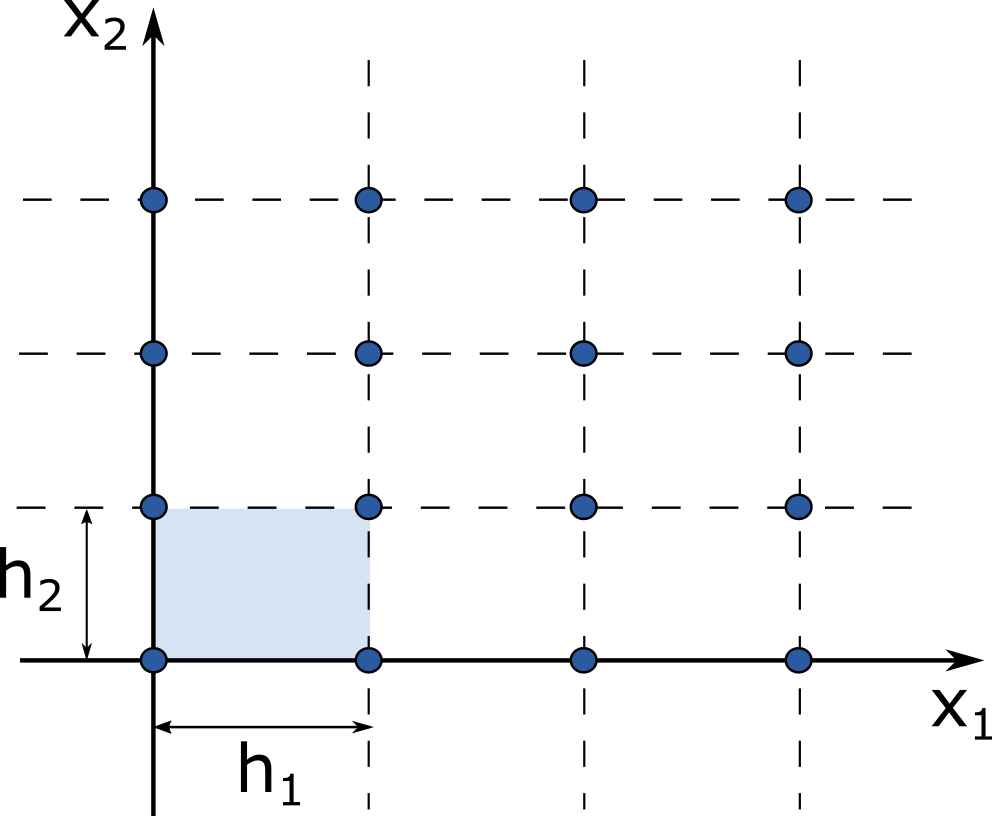}
   \caption{A design of experiments $D_H$ for $\iD = 2$.}
   \label{fig:grid_design}
\end{figure}

We observe the random process at the infinite grid $D_H = \{ x_{\vecK} = H \vecK, k \in \mathbb{Z}^{\iD}\}$.
$H$ is a diagonal matrix with elements at the diagonal $\mathrm{diag}\{h_1, \ldots, h_{\iD}\}$.
An example of such two dimensional design of experiments is given at Figure~\ref{fig:grid_design}.

We investigate interpolation error of $f(\vecX)$
using the best regression model $\tilde{f}(\vecX)$.
In Gaussian case this model depends linearly on observations
\[
\tilde{f}(\vecX) = \sum_{\vecK \in \mathbb{Z}^{\iD}} K(\vecX - \vecX_\vecK) f(\vecX_k),
\]
where $K(\vecX)$ is a kernel function obtained as a solution of Kolmogorov-Wiener--Hopf equations.

For a set 
$\Omega_H = \prod_{i = 1}^\iD [0, h_i]$
we are interested in evaluation of the integral of the expectation of 
squared differences between true value of a random process and his interpolation:
\[
\sigma^2(\tilde{f}, F) = \frac{1}{\Omega_H} \int_{\Omega_H} \mathbb{E} \left[\tilde{f}(\vecX) - f(\vecX) \right]^2 d \vecX.
\] 

The following theorem holds:
\begin{theorem}
\label{th:error_u}
For a random process $f(\vecX)$ with a spectral density $F(\vecO)$, observed at $D_H$
the interpolation error has the form:
\begin{equation*}
\label{eq:interpolation_error}
\sigma^2(\tilde{f}, F) = \int_{\Omega_H} F(\vecO) \left[\left(1 - \hat{K}(\vecO)\right)^2 + \sum_{\vecK \ne 0 } \hat{K}^2 \left(\vecO + H^{-1} \vecK \right)\right] d\vecO,
\end{equation*}
where $\hat{K}(\vecO)$ is a Fourier transform of $K(\vecO)$.
Moreover, $\hat{K}(\vecO)$ has the form
\[
\hat{K}(\vecO) = F(\vecO) / \sum_{\vecK} F \left(\vecO + H^{-1} \vecK \right).
\]
\end{theorem}

Often the true spectral density is unknown.
So, we are interested in the minimax interpolation error:
\begin{equation}
\label{eq:minimax_error}
R^H(L) = \inf_{\tilde{f}} \sup_{F \in \mathcal{F}(L)} \sigma^2(\tilde{f}, F),
\end{equation}
where $\mathcal{F}(L)$ defines a set of spectral densities that correspond to smooth enough Gaussian processes:
\[
\sup_{F \in \mathcal{F}(L)} \mathbb{E} \| f^{(1)}(\vecX) \|^2 \leq L, 
\]
$f^{(1)}(\vecX)$ is a vector of first partial derivatives of Gaussian process with a spectral density $F(\vecO)$.

The following theorem holds:
\begin{theorem}
\label{th:minimax_error}
The minimax interpolation error $R^H(L)$ from~\eqref{eq:minimax_error} has the form:
\[
R^H(L) = \frac{L}{2 \pi^2} \max_{i = \overline{1, \iD}} h_i^2.
\]
\end{theorem}

Given the theorems above we can get the interpolation error for certain covariance functions: exponential and squared exponential covariance functions.

\begin{corollary} \label{col:exponential_error}
For Gaussian process at $\bbR$ with the exponential spectral density of the form $F_{\theta}(\omega) = \frac{\theta}{\theta^2 + \omega^2}$
the interpolation error~\eqref{eq:interpolation_error} for the best interpolation has the form:
\begin{equation*}
    \sigma_h^2(\tilde{f}, F_{\theta}) \approx \frac{2}{3} \pi^2 \theta h + O((\theta h)^2),
\end{equation*}
for $\theta h \rightarrow 0$.
\end{corollary}

\begin{corollary} \label{col:squared_exponential_error}
For Gaussian process at $\bbR$ with the squared exponential spectral density of the form $F_{\theta}(\omega) = \frac{1}{\sqrt{\theta}} \exp \left(-\frac{\omega^2}{2 \theta}\right)$
the interpolation error~\eqref{eq:interpolation_error} for the best interpolation has the form:
\begin{align*}
    &\frac{4}{3} h \sqrt{\theta} \exp \left(-\frac{1}{8 h^2 \theta}\right)
        \leq \sigma^2_{h}(\tilde{f}, F_{\theta}) \leq 7 h \sqrt{\theta} \exp \left(-\frac{1}{8 h^2 \theta} \right)
\end{align*}
for $\theta h^2 \rightarrow 0$.
\end{corollary}

So, the minimax error decreases as $h^2$ for $h = \max_{i \in \{1, \iD\}} h_i$,
while for some covariance functions it can decrease exponential with respect to $h$, or decrease linearly with $h$ for a non-smooth Gaussian process.

\section{Interpolation error for misspecified case}
\label{sec:wrong}


In practice we use a model of the true Gaussian process.
Let us consider a Gaussian process with the true spectral density $F(\vecO)$, while for estimation we use a Gaussian process
with the spectral density $F_{\theta}(\vecO)$.
The problem is to estimate the interpolation error for misspecified spectral density used for computation of the final approximation.

We again consider the infinite grid design of experiments $D_H$ 
and sample of values $\{ f(\vecX_{\vecK}) \}$ at $D_H$ of a realization of a Gaussian process with the spectral density $F(\vecO)$.

The best interpolation has the form:
\begin{align*}
\tilde{f}_{\theta}(\vecX) &= H \sum_{\vecK \in \mathbb{Z}^{\iD}} K_{\theta}(\vecX - \vecX_{\vecK})  f(\vecX_{\vecK}).
\end{align*}
We obtain the kernel $K(\cdot)$ by minimization of the mean squared error assuming that the true spectral density is $F(\vecO)$.
We obtain the kernel $K_{\theta}(\cdot)$ in a similar way, 
but using the true spectral density  $F_{\theta}(\vecO)$.

Our goal is to estimate the interpolation error
\[
\sigma^2_{H}(\tilde{f}_{\theta}, F) = \frac{1}{\Omega_H} \int_{\Omega_H} \mathbb{E} \left[\tilde{f}_{\theta}(\vecX) - f(\vecX) \right]^2 d\vecX.
\]

\begin{theorem}
\label{th:mis}
The interpolation error for the true spectral density $F(\vecO)$, 
if we used the spectral density $F_{\theta}(\vecO)$ for construction of the regression model given observations at $H^{-1} \vecK, \vecK \in \mathbb{Z}^{\iD}$ has the form:
\begin{equation*}
\label{eq:interpolation_error}
\sigma^2_{H}(\tilde{f}_{\theta}, F) = \int_{\Omega_H} F(\vecO) \left[\left(1 - \hat{K}_{\theta}(\vecO)\right)^2 + \sum_{\vecK \ne 0 } \hat{K}_{\theta}^2 \left(\vecO + H^{-1} \vecK \right)\right] d\vecO.
\end{equation*}
\end{theorem}

So, given spectral densities $F(\vecO)$ and $F_{\theta}(\vecO)$ 
we can get the target interpolation error by analytical integration of \eqref{eq:interpolation_error} or numerical estimation.

Note, that we can get the result of Theorem~\ref{th:mis} in the form:   
\[
\sigma^2_{H}(\tilde{f}_{\theta}, F) = \int_{\Omega_H} F(\vecO)
\left\{ \frac{\sum_{\vecK \ne 0} F(\vecO + H^{-1} \vecK)}{\sum_{\vecS} F(\vecO + H^{-1} \vecS)} \right\} d\vecO
\]

As a set of coefficients $K(\vecX - \vecX_{\vecK})$ minimizes the interpolation error, it holds that $\sigma^2_{H}(\tilde{f}_{\theta}, F) \geq \sigma^2_{H}(\tilde{f}, F)$.
Now we are ready for analysis of difference of the interpolation error in the cases of misspecified and correctly specified models.

\subsection{Interpolation error for misspecified Matern spectral density $\nu = \frac12$}
\label{sec:ornstein}

We consider the interpolation error for the misspecified case 
for Matern with $\nu = \frac12$ spectral density.
For $\bbR$ and a stationary Gaussian process with Matern $\nu = \frac12$ covariance function $R_{\theta}(x)$:
\[
R_{\theta}(x) = \sqrt{\frac{\pi}{2}} \exp \left(-\theta \left\|x\right\| \right).
\]
An alternative name for this covariance function is the exponential covariance function~\cite{gpy2014}.
The spectral density that corresponds to this covariance function is the following:
\[
F_{\theta}(\omega) = \frac{\theta}{\theta^2 + \omega^2}.
\]

To construct an interpolation we use a misspecified spectral density  $F_{\wtheta}, \wtheta \ne \theta$.

\begin{corollary} \label{col:exponential_error_misspec}
We observe a realization of Gaussian process at $D_h \subset \bbR$ 
with the true exponential spectral density of the form $F_{\theta}(\omega) = \frac{\theta}{\theta^2 + \omega^2}$.
Then the interpolation error~\eqref{eq:interpolation_error} for $\tilde{f}_{\wtheta}$ constructed with the assumption that the true spectral density is $F_{\wtheta}(\omega)$ has the form:
\begin{equation*}
    \sigma_h^2(\tilde{f}_{\wtheta}, F_{\theta}) \approx \frac{2}{3} \pi^2 \theta h + O((\theta h)^2),\,\, 
    \theta h \rightarrow 0 \,.
\end{equation*}
\end{corollary}

So, for small $h$ the interpolation error doesn't depend on the coefficient $\wtheta$.
It is obvious that generally the misspecified case lacks this nice property.

\section{Computational experiments}
\label{sec:experiments}

We obtained theoretical results in section \ref{sec:prior} under assumption that the realization of a Gaussian process is known at an infinite grid. 
It is impossible to fulfill such assumption in practice. 
While, we expect that results will be the same for a large enough finite sample and an infinite grid of points,
we should validate theoretical for a finite sample.

For this purpose Gaussian process realization with Matern covariance function with $\nu = \frac12$ is taken as objective function. 
Matern covariance function with $\nu = \frac12$ specifies as follows: $R_\theta(x) = \sqrt[]{\frac{\pi}{2}}\exp(-\theta||x||)$. It is a special case of exponential covariance function, so, results of section \ref{sec:prior} are applicable for it.

\subsection{Workflow of computational experiments}

In this subsection we provide technical details on computational experiments.
For experiments we used Gaussian process regression realization from \cite{gpy2014} library. 
We provide the code used to generate results in the article at \cite{github} github page. 
For the sake of clarity and faster convergence of the empirical interpolation error to the true one we consider one-dimensional grids of points, while our theoretical results are valid for the multivariate case.

There are three steps in computational experiment dedicated to obtaining of the interpolation error: create a realization of a Gaussian process; use a regression model with an alternative covariance function on the base of a training sample; estimate the interpolation error given the constructed regression model and a test sample.

We create realizations of a Gaussian process using the following steps:

\begin{enumerate}
\item  Generate a grid of points $X$ of size $\sS$ from interval $[0, 1]$. Step of the grid $h$ is inversely proportional to the current sample size $\sS$.
\item  Select the parameter of covariance function $\theta$.
\item  Evaluate the sample covariance matrix $K$ for points $X$. We use Matern covariance function with parameter $\theta$ as a covariance function and white noise with variance $10^{-8}$.
\item  Evaluate the Cholesky decomposition $L$ of the obtained covariance matrix $K$.
\item  Generate a vector of i.i.d. random variables $\vecY_0$ from the standard normal distribution of size $\sS$.
\item  Obtain multivariate normal distribution by multiplying the Cholesky decomposition $L$ by $\vecY_0$  at the previous step.
\end{enumerate}

\begin{figure}
\subfigure[Smoothness 10]{\label{fig:gen10}\includegraphics[width=50mm]{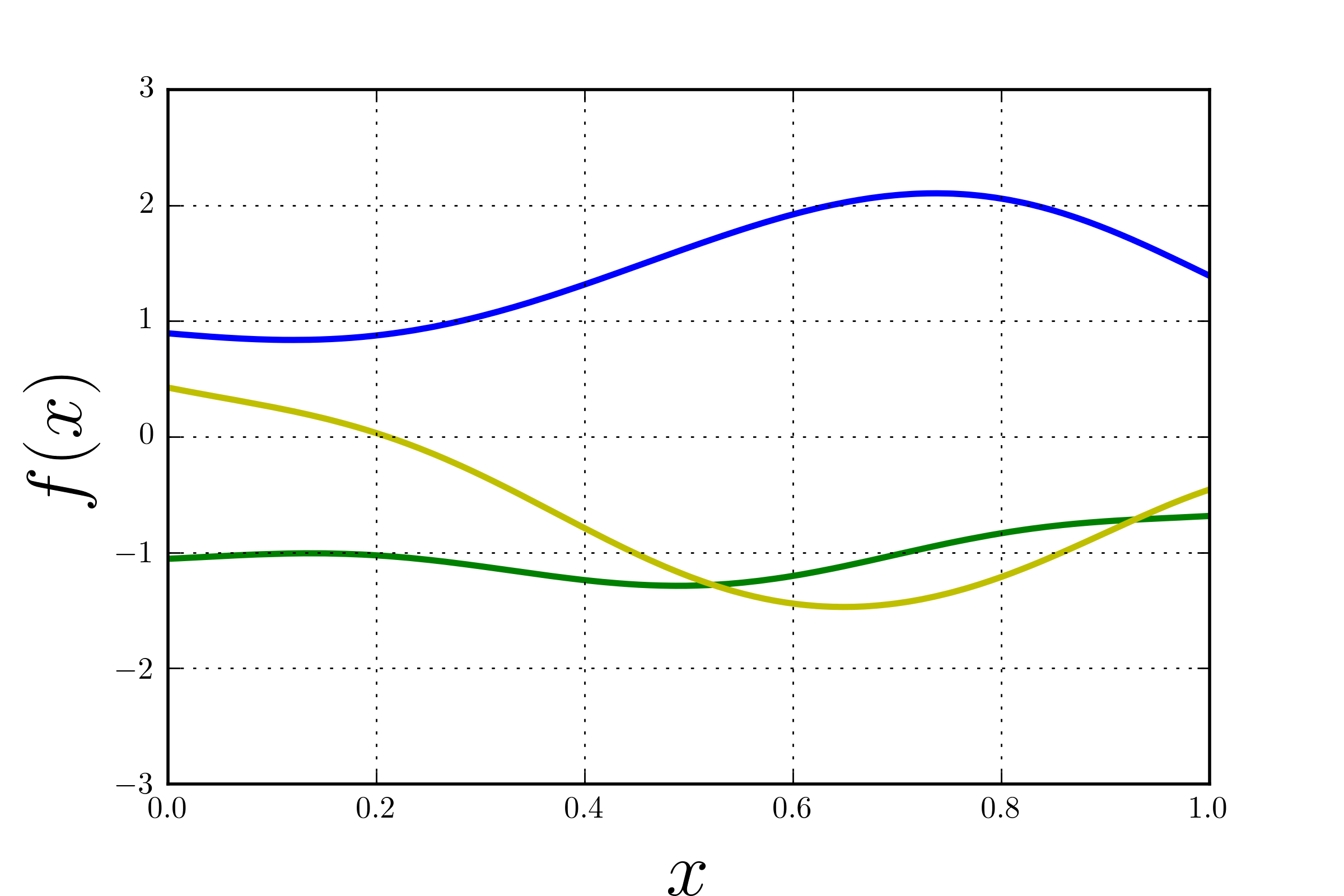}}
\subfigure[Smoothness 100]{\label{fig:gen100}\includegraphics[width=50mm]{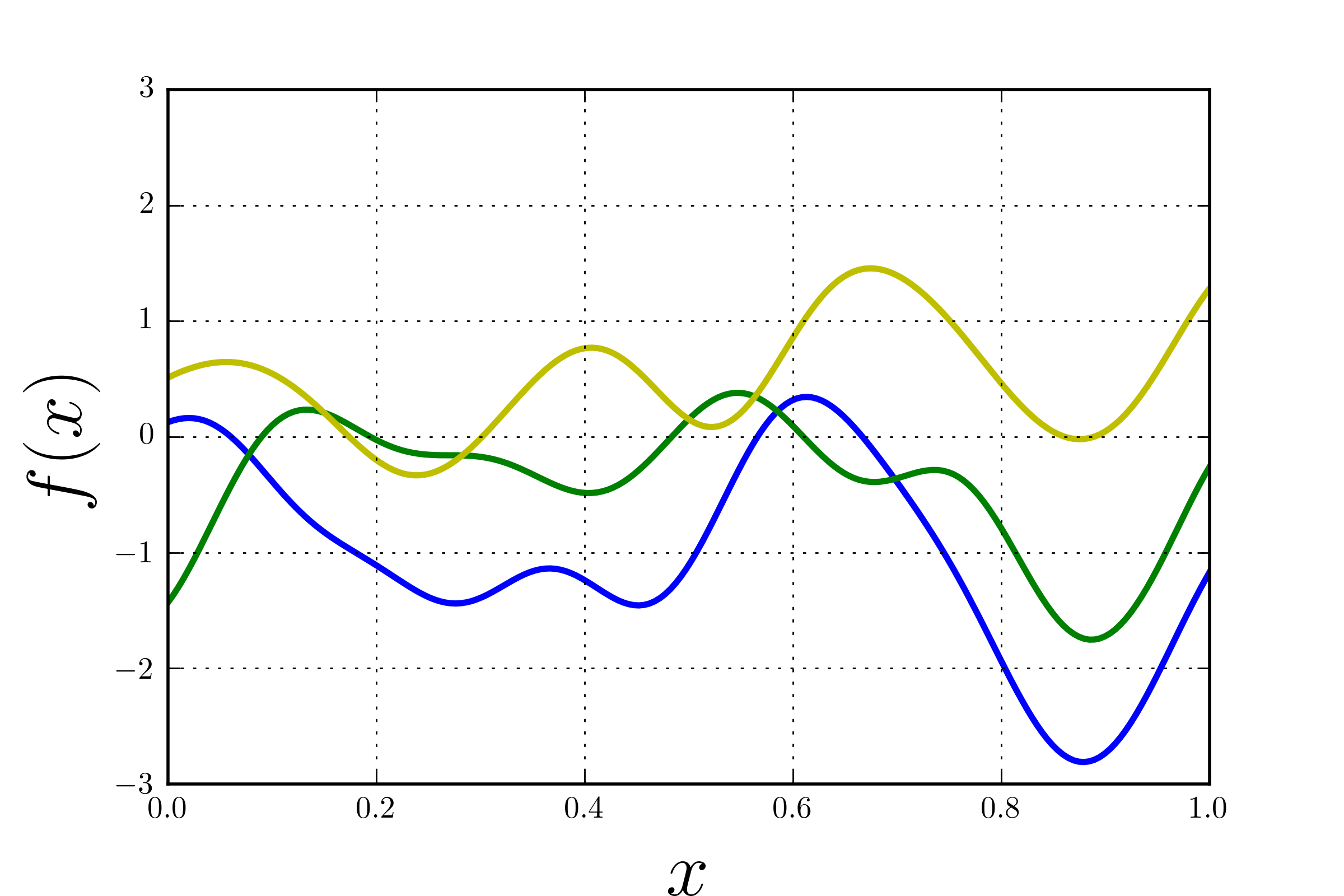}}
\subfigure[Smoothness 500]{\label{fig:gen500}\includegraphics[width=50mm]{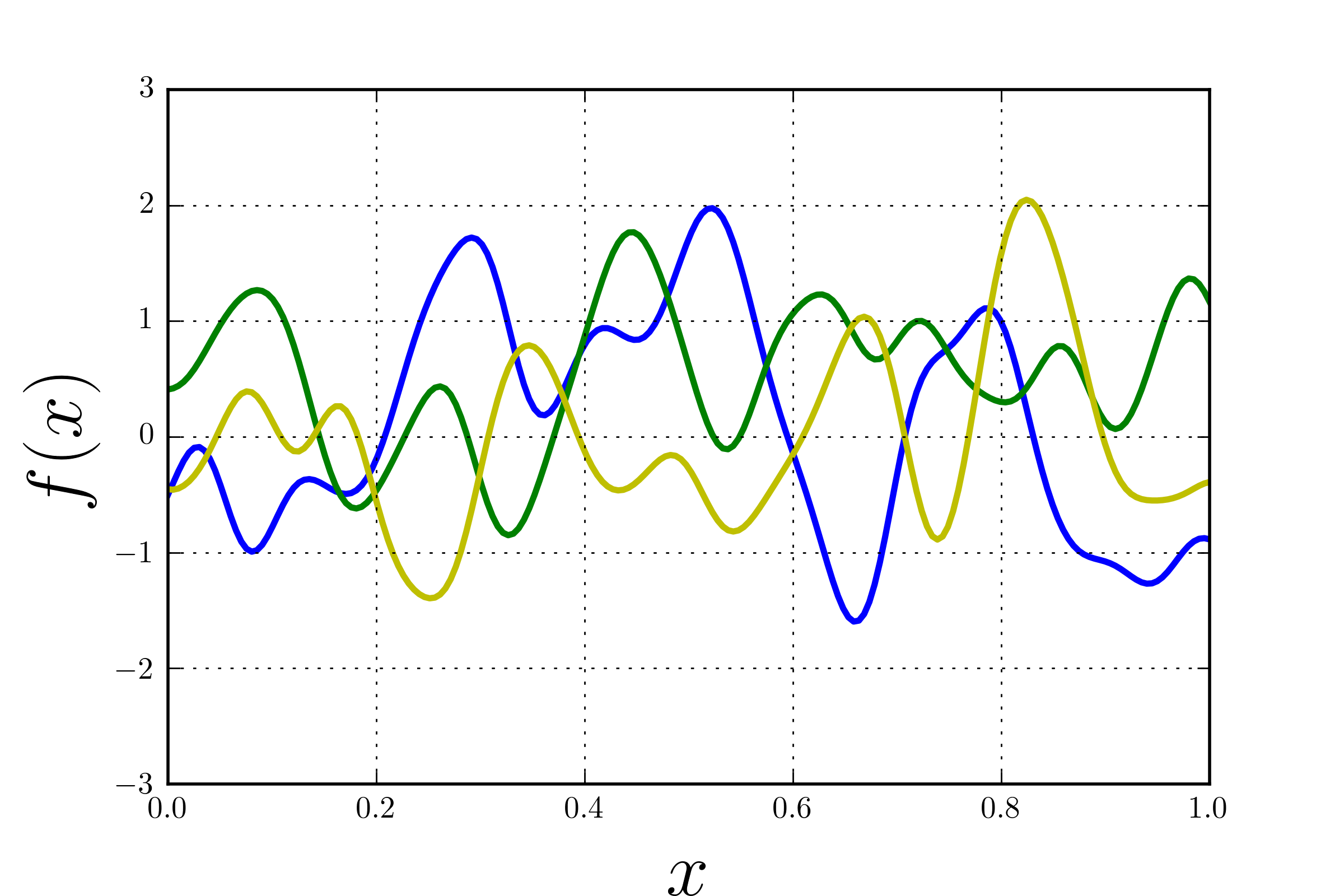}}
  \caption{Realizations of Gaussian process with the squared exponential covariance function for different smoothness values}
  \label{fig:data_gen}
\end{figure}

As the results of this procedure we get $\vecY = C \vecY_0 \sim \mathcal{N}(\mathbf{0}, K)$. 
Few examples of generated realizations are at Figures~\ref{fig:data_gen}.

Next algorithm is used for construction of the regression model and evaluation of the interpolation error of this model:

\begin{enumerate}
\item  Select the parameter $\wtheta$ --- an assumption about the true parameter of the covariance function.
\item  Create a covariance function with chosen parameters.
\item  Specify regression model on grid, according to the kernel with chosen parameters. 
\item  Estimate the interpolation error for the selected sample size as $\sum_{i = 1}^{\sS_{\mathrm{test}}}\left(\hat y_i - y_i \right)^2$, where $\hat y_i$ and $y_i$ are respectively the predicted and the true value at a test point $\vecX_i$ using a test sample with a more dense grid.
\end{enumerate}

\subsection{Simulations}

\begin{figure}[h]
\centering    
\subfigure[A single run of the algorithm]{\label{fig:choletskiy1}\includegraphics[width=0.45\linewidth]{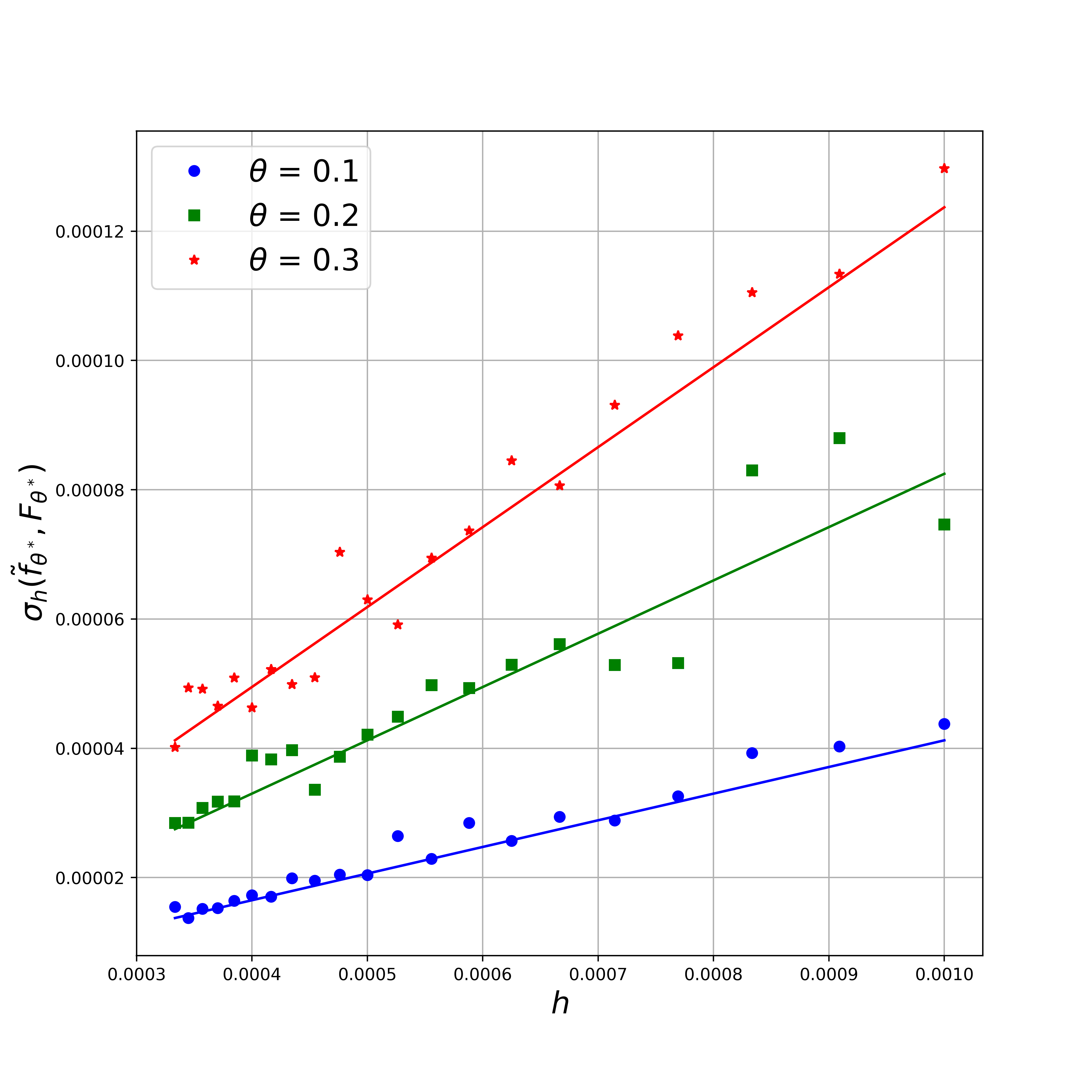}}
\subfigure[Averaging of ten runs]{\label{fig:choletskiy2}\includegraphics[width=0.45\linewidth]{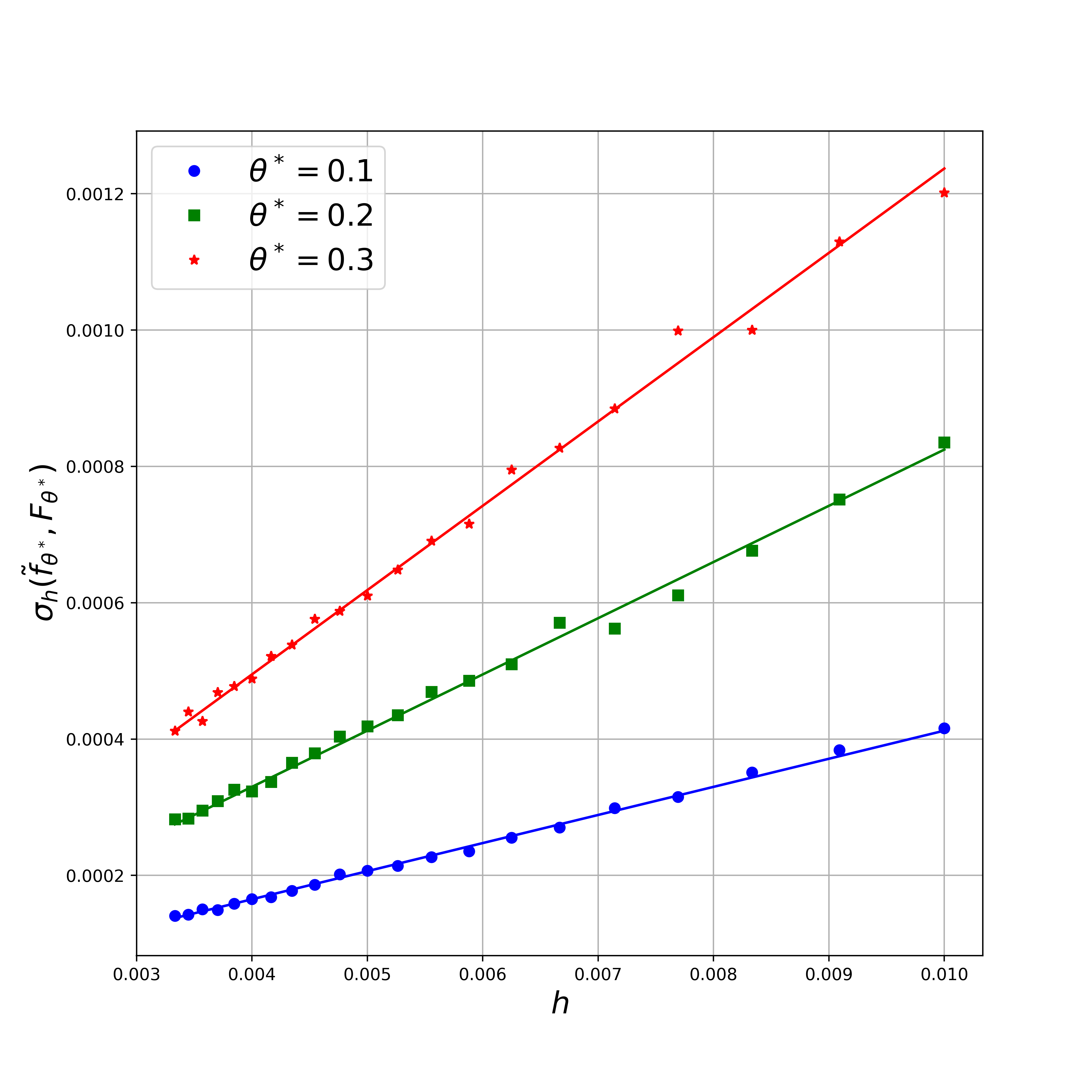}}
\caption{The interpolation errors obtained in the experiment. The solid lines indicate theoretical results.}\label{choletski}
\end{figure}

We start with the following setup: $\theta = \wtheta$ for values $\theta = [0.1, 0.2, 0.3]$. 
Obtained results are at Figure~\ref{choletski}.
For each of them the "basic algorithm" had been launched. 
After that the experiment we average the result of 10 realizations. The obtained interpolation errors fit the straight line even better.
So, the assumption of the infinite grid does not affect the interpolation error: the theoretical and the practical results are consistent.

We continue with the misspecified problem:
for a fixed $\theta = 0.1$, three different $\wtheta = [0.1, 1, 10]$ are used during model construction. 
For each of the values at each sample size, we run the "basic algorithm" 20 times and average the results.
We see  that results are similar for different used values of $h$. While the obtained results slightly differ, Theorem~\ref{th:error_u} gives almost perfect approximation in this case.

\begin{figure}[h]
\begin{center}
\begin{tabular}{lccc}
\hline
h & $\theta$ = 0.1 & $\theta$ = 0.1 & $\theta$ = 0.1 \\
 & $\theta'$ = 0.1 & $\theta'$ = 1 & $\theta'$ = 10\\
\hline
0.01 & 0.425 & 0.413 & 0.424\\
0.004 & 0.166 & 0.168 & 0.162\\
0.0025 & 0.102 & 0.105 & 0.104\\
\hline
\end{tabular}
\caption{Obtained values of the interpolation error $\sigma^2_h(\hat{f}_\theta, F_\theta) \cdot 10^3$ for misspecified case averaged over 20 realizations}
\end{center}
\end{figure}

\begin{figure}[h]
\centering    
\subfigure[The exponential covariance function]{\label{fig:exponential}\includegraphics[width=0.31\linewidth]{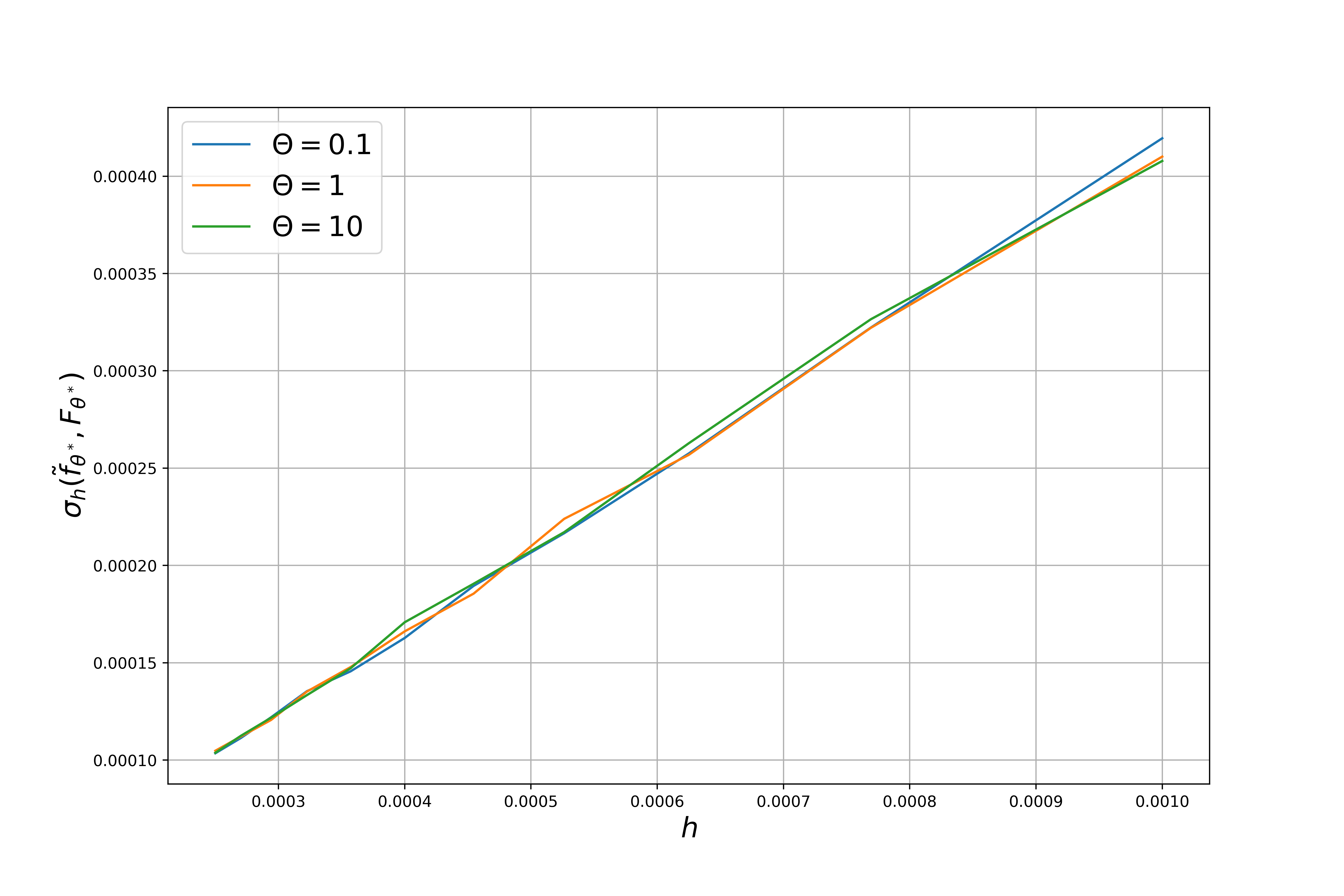}}
\subfigure[$Matern_{3/2}$ covariance function]{\label{fig:matern_32}\includegraphics[width=0.31\linewidth]{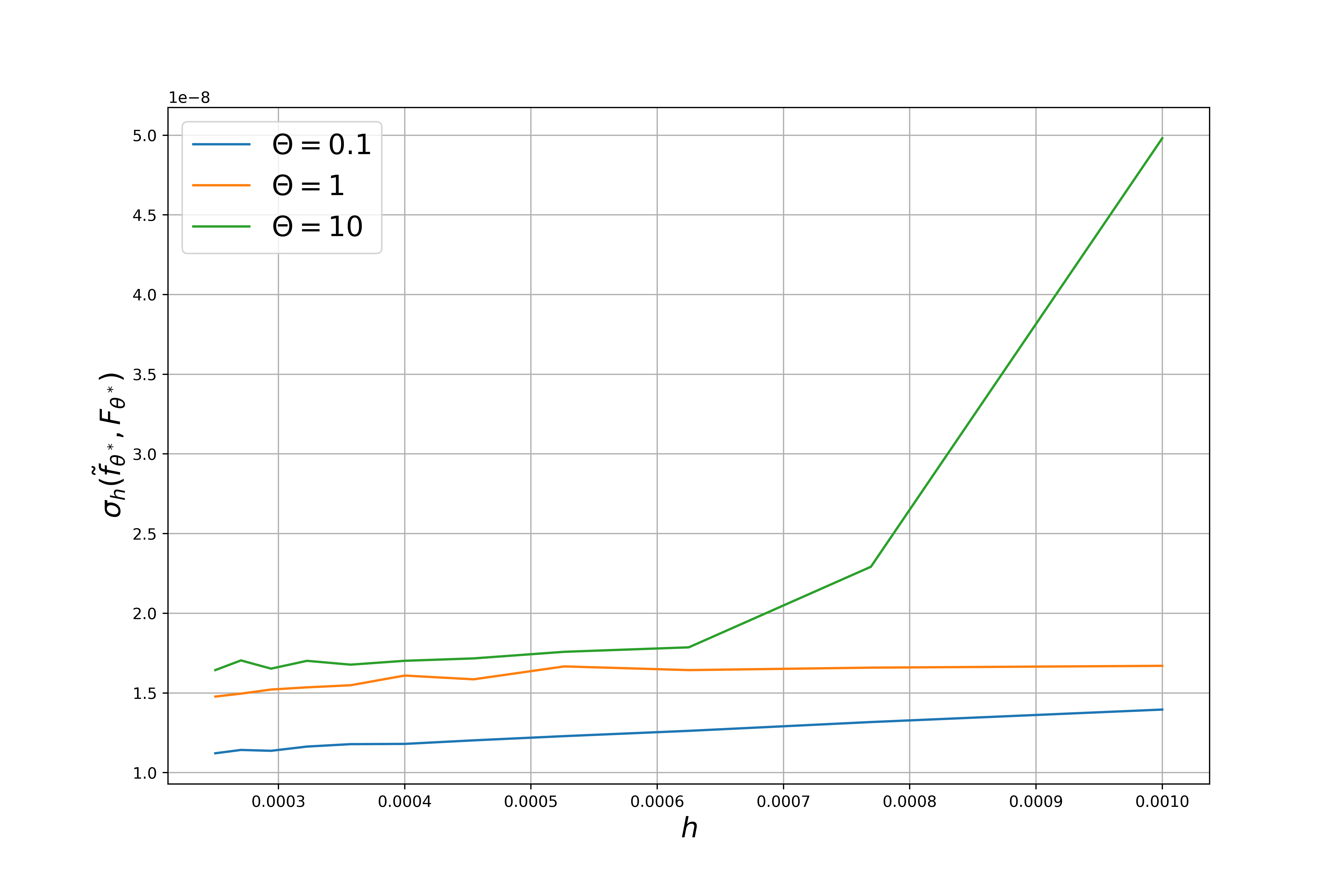}}
\subfigure[The squared exponential covariance function]{\label{fig:matern_RBF}\includegraphics[width=0.31\linewidth]{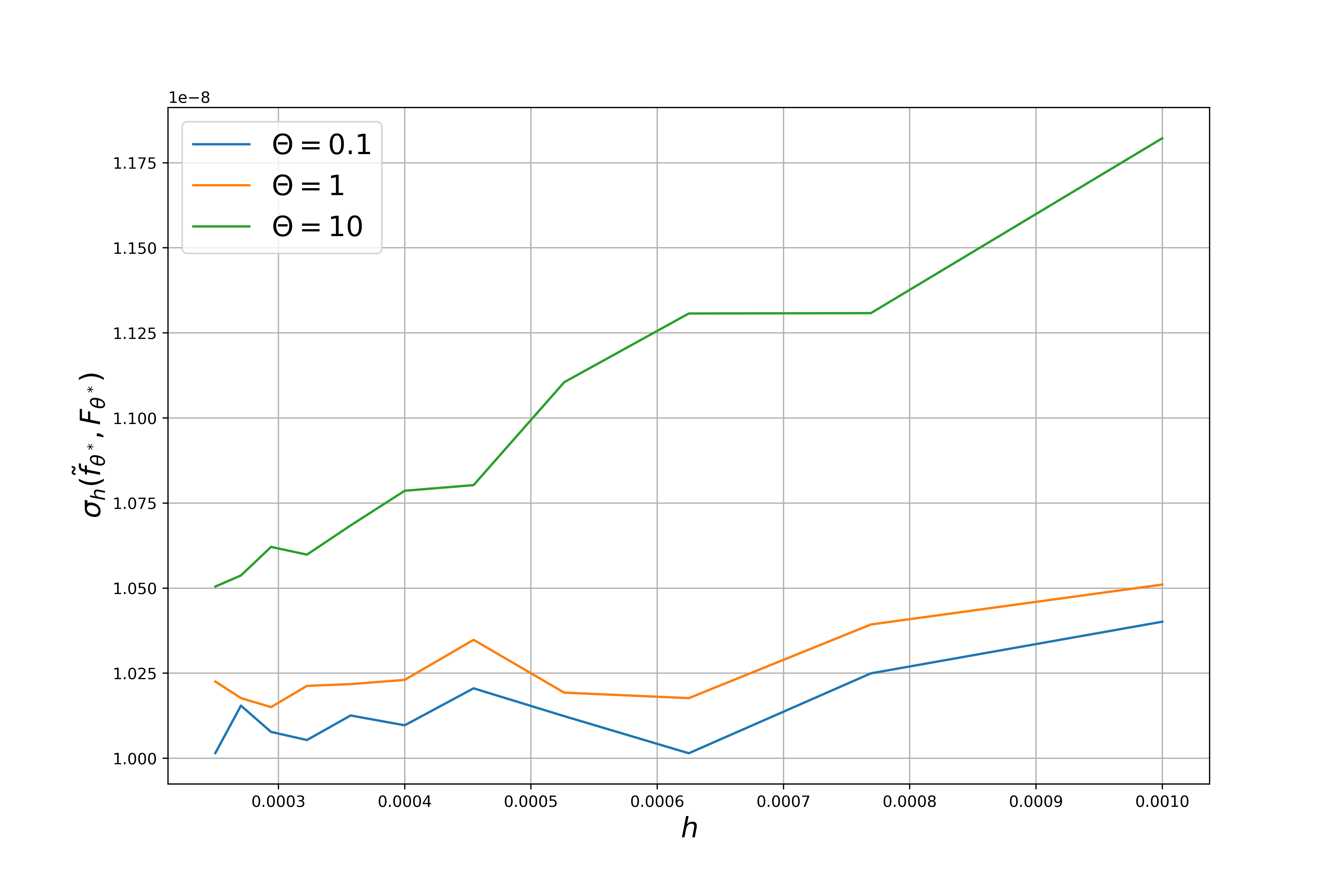}}
\caption{The interpolation errors obtained in the experiment for different classes of covariance functions averaged over $20$ realizations. The true parameter is $1$ in all experiments.}
\label{fig:various_cf}
\end{figure}

For the exponential, $Matern_{3/2}$ and RBF kernels we examine the interpolation error for the misspecified model, that lies in the same parametric class of models.
The obtained interpolation errors are at Figure~\ref{fig:various_cf}.
We also run the Wilcoxon difference test~\cite{wilcoxon1945individual} to test if the results are different.
For the exponential kernel we get $p-\mathrm{value} > 0.92$, while for $Matern_{3/2}$ and the squared exponential kernel $p-\mathrm{value} < 10^{-2}$. 
So, for the exponential kernel we can't reject the hypothesis that the interpolation errors are the same, no matter what value of the covariance function parameter is used, while for two other covariance functions results suggest that the hypothesis that the interpolation errors are the same seems to be wrong.
Therefore, for a non-exponential kernel the interpolation error depends on the value of the parameter used for the construction of the regression model.

\begin{figure}[!h]
\centering
\includegraphics[width = 0.7\linewidth]{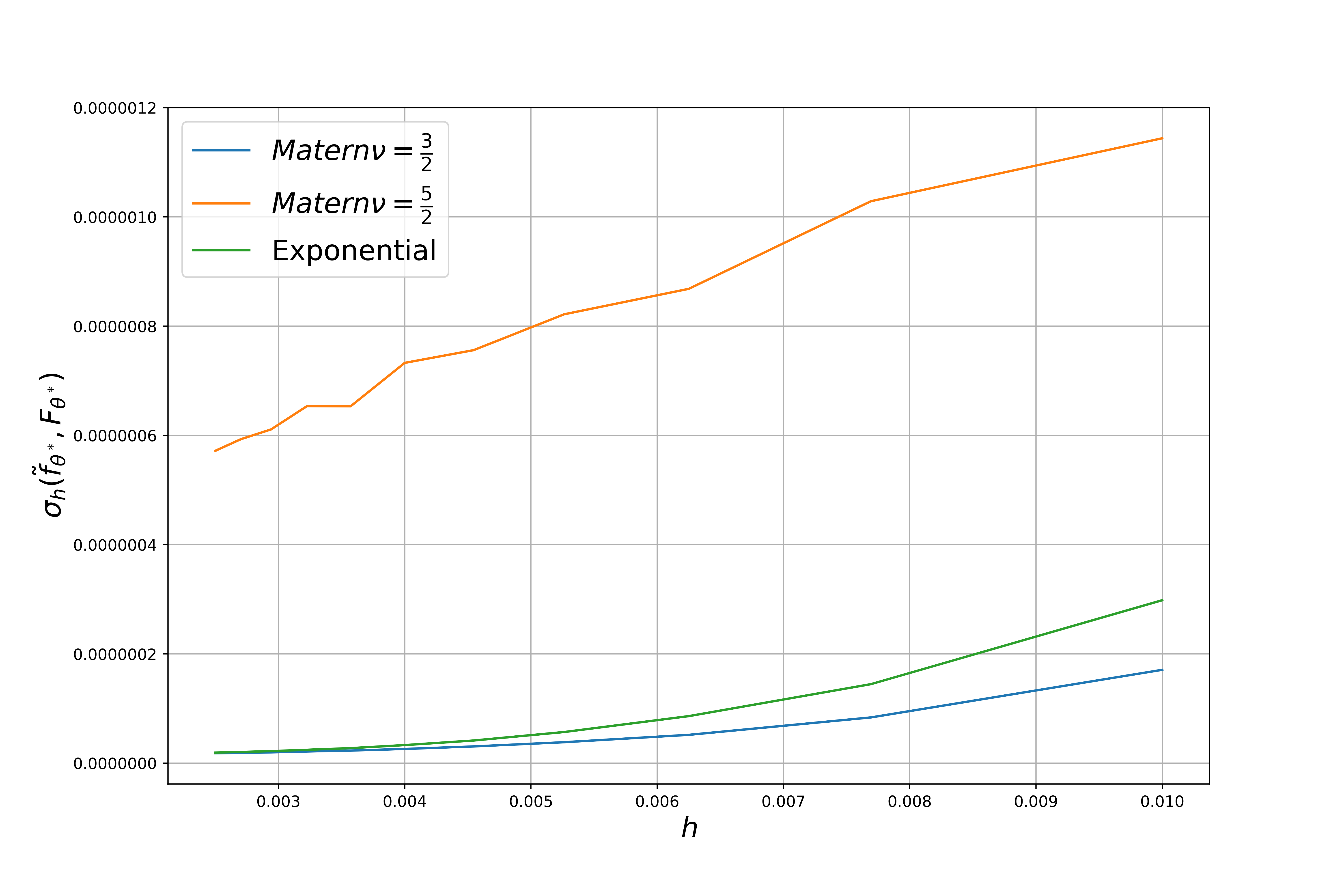}
\caption{The Interpolation errors obtained in the experiment for misspecified covariance functions averaged over 20 realizations.}
\label{fig:wrong_family}
\end{figure}

In the last experiment, we investigate the interpolation error for the case when the wrong parametric class of functions is used.
In particular, the parameter of the model is $\theta = 1$;
when the model $Matern_{3/2}$ is true, models with covariance functions $Matern_{3/2}$, $Matern_{5/2}$ and exponential give the interpolation errors provided at Figure~\ref{fig:wrong_family}.
We see that the error estimation derived in Corollary~\ref{col:exponential_error} is not applicable in the case of using different classes of the true and used function.

\section{Conclusions}

This article presents the interpolation error for misspecified regression model of Gaussian process regression.
The obtained result can be used for analysis of effect of the model misspecification on the quality of obtained regression model.
For example, for Matern covariance function the interpolation error doesn't depend on used value of parameter.
This effect holds for numerical experiments for non-grid finite training samples.

\acks{The results presented in the sections \ref{sec:wrong}, \ref{sec:ornstein} were supported solely by the Russian Science Foundation grant (project 14-50-00150).}

\bibliography{main}

\appendix

\section{Proofs of the presented results}\label{apd:first}

The proof of Theorem~\ref{th:mis}:
\begin{proof}
The difference between this problem and the problem given in Theorem~\ref{th:error_u} is on different set of coefficient $K_{\theta}(\vecX - \vecX_{\vecK})$ identified by the spectral density.
Consequently, we are able to get the proof of Theorem~\ref{th:error_u} using the results given in~\cite{zaytsev2017minimax}.

It is easy to see that
\begin{align*}
&\bbE [f(\vecX) - \tilde{f}_{\theta}(\vecX) ]^2 = \int_{\bbR^{\iD}} F(\vecO) \left|1 - |H| \sum_{\vecK \in \bbZ^{\iD}} K_{\theta}(\vecX - \vecX_{\vecK}) \exp(-2\pi \mathrm{i} \vecO^{T} (\vecX_{\vecK} - \vecX)) \right|^2 d\vecO = \\
&= \int_{\bbR^{\iD}} F(\vecO) \left|1 - |H| \sum_{\vecK \in \bbZ^{\iD}} \left( \int_{\bbR^{\iD}} \hat{K}_{\theta}(\vecU) \exp(-2\pi \mathrm{i} \vecU^T (\vecX - \vecX_{\vecK})) d\vecU \right) \exp(-2\pi \mathrm{i} \vecO^{T} (\vecX_{\vecK} - \vecX)) \right|^2 d\vecO,
\end{align*}
where $\hat{K}_{\theta}(\vecU)$ is the Fourier transform of $K_{\theta}(\vecX)$.
As Poisson summation formula suggests:
\[
\sum_{\vecK \in \bbZ^{\iD}} \exp(2 \pi \mathrm{i} \vecK^T \vecO) = \sum_{\vecK \in \bbZ^{\iD}} \delta(\vecO + \vecK),
\]
where $\delta(\vecO)$ is the Dirac delta function,
then
\begin{align*}
&\bbE [f(\vecX) - \tilde{f}_{\theta}(\vecX) ]^2 = 
\int_{\bbR^{\iD}} F(\vecO) \left|1 - |H| \sum_{\vecK \in \bbZ^{\iD}} \int_{\bbR^{\iD}} \hat{K}_{\theta}(\vecU) \exp(2\pi \mathrm{i} (\vecO - \vecU)^T \vecX) \delta(\vecU - \vecO + H^{-1} \vecK)d\vecU \right|^2 d\vecO = \\
&= \int_{\bbR^{\iD}} F(\vecO) \left| 1 - \sum_{\vecK \in \bbR^{\iD}} \hat{K}_{\theta} 
(\vecO - H^{-1} \vecK) \exp(2 \pi \mathrm{i} H^{-1} \vecX^T \vecK) \right|^2 d\vecO.
\end{align*}

Taking into account orthogonality of the system of functions $\exp(2 \pi \mathrm{i} H^{-1} \vecX^T \vecK)$ on $\vecX \in [0, h_1] \times \ldots \times [0, h_{\iD}]$ we integrate the equality to get the interpolation error
\[
\sigma^2_H(\tilde{f}_{\theta}, F) = 
\int_{\bbR^{\iD}} F(\vecO) \left|[1 - \hat{K}_{\theta} (\vecO)]^2 + \sum_{\vecK \in \bbZ^{\iD} \setminus \{\mathbf{0}\}} \hat{K}_{\theta} ^2 (\vecO + H^{-1} \vecK) \right|^2 d\vecO.
\]

To get $\hat{K}_{\theta} (\vecO)$ that minimizes the interpolation error we rewrite the equation above using $F_{\theta}$ instead of $F$: 
\[
\sigma^2_H(\tilde{f}_{\theta} , F_{\theta}) = 
\int_{\bbR^{\iD}} \left|[1 - \hat{K}_{\theta}(\vecO)]^2 F_{\theta}(\vecO) + \hat{K}_{\theta}(\vecO)^2 \sum_{\vecK \in \bbZ^{\iD} \setminus \{\mathbf{0}\}} \hat{F}_{\theta} (\vecO + H^{-1} \vecK) \right|^2 d\vecO.
\]
To minimize this error we solve this quadratic optimization problem for each $\vecO$ and get:
\[
\hat{K}_{\theta}(\vecO) = \frac{\hat{F}_{\theta} (\vecO)}
                      {\sum_{\vecK \in \bbZ^{\iD}} \hat{F}_{\theta} (\vecO + H^{-1} \vecK)}.
\]
Then 
\begin{equation}
\sigma^2_H(\tilde{f}_{\theta}, F) = \int_{\bbR^{\iD}} F(\vecO) \frac{\sum_{\vecK \in \bbZ^{\iD} \setminus \{\mathbf{0}\}} \hat{F}_{\theta} (\vecO + H^{-1} \vecK)}
                      {\sum_{\vecK \in \bbZ^{\iD}} \hat{F}_{\theta} (\vecO + H^{-1} \vecK)} d\vecO.
\end{equation}
\end{proof}

\noindent The proof of Corollary~\ref{col:exponential_error_misspec}:
\begin{proof}
Our goal is to evaluate
\[
\sigma_h^2 \left(\tilde{f}_{\wtheta}, F_{\theta}\right) = \int_{-\infty}^{\infty} F_{\theta}(\omega) \frac{\sum_{k \ne 0 } F_{\wtheta}(\omega + \frac{k}{h})}{\sum_k F_{\wtheta}(\omega + \frac{k}{h})} d\omega.
\]

It holds that
\begin{align*}
&\sum_{k = -\infty}^{\infty} F_{\theta} \left(\omega + \frac{k}{h} \right) = \sum_{k = -\infty}^{\infty} \frac{\theta}{(\omega + \frac{k}{h})^2 + \theta^2} = 
h \sum_{k = -\infty}^{\infty} \frac{h \theta}{(h \omega + k)^2 + h^2 \theta^2} = \\
&=\pi h \coth(\pi \theta h) \frac{1}{1 + \sin^2(\pi h \omega) (\coth^2(\pi \theta h) - 1)}.
\end{align*}

Then
\[
\int_{-\infty}^{\infty} F_{\theta}(\omega) \frac{\sum_{k \ne 0 } F_{\wtheta}(\omega + \frac{k}{h})}{\sum_k F_{\wtheta}(\omega + \frac{k}{h})} d\omega = \int_{-\infty}^{\infty} \frac{\theta}{\theta^2 + \omega^2} \left(1 - \frac{\wtheta}{\wtheta^2 + \omega^2} \frac{1 + \sin^2(\pi h \omega) (\coth^2(\pi \wtheta h) - 1)}{\pi h \coth(\pi \wtheta h)} \right) d\omega.
\]

For three integrals presented above it holds that:
\[
\int_{-\infty}^{\infty} \frac{\theta}{\theta^2 + \omega^2} d\omega = \pi.
\]
Moreover,
\[
\int_{-\infty}^{\infty} \frac{\theta}{(\theta^2 + \omega^2)} \frac{\wtheta}{(\wtheta^2 + \omega^2)} d\omega = \frac{\pi}{\theta + \wtheta}.
\]
Finally,
\[
\int_{-\infty}^{\infty} \frac{\theta}{(\theta^2 + \omega^2)} \frac{\wtheta}{(\wtheta^2 + \omega^2)} \sin^2(\pi \omega h) d\omega = \frac{\pi}{2(\theta + \wtheta)} \left(1 - \frac{\theta \exp(-2 \pi h \wtheta) - \wtheta \exp(-2 \pi h \theta)}{\theta - \wtheta} \right).
\]

Consequently
\begin{align*}
&\sigma_h^2 \left(\tilde{f}_{\wtheta}, F_{\theta}\right) = 
\pi - \frac{\pi}{2 \pi (\theta + \wtheta) h \coth(\pi \wtheta h)} + \\
&+ \frac{\pi}{2 \pi (\theta + \wtheta) h \coth(\pi \wtheta h)} \left(1 - \frac{\theta \exp(-2 \pi h \wtheta) - \wtheta \exp(-2 \pi h \theta)}{\theta - \wtheta} \right) \frac{\coth^2(\pi \wtheta h) - 1}{\pi h \coth(\pi \wtheta h)}.
\end{align*}

We are interesting in case $h \rightarrow 0$.
In this case we can use Taylor decomposition to evaluate the final result.
We get it in the following form:
\[
\sigma_h^2 \left(\tilde{f}_{\wtheta}, F_{\theta}\right) = \frac{2 \pi^2}{3} \theta h + O((\theta h)^2) + O((\wtheta h)^2).
\]
\end{proof}

\end{document}

%% file: mydef.tex
\newcommand\vecK{\mathbf{k}}

\newcommand\vecS{\mathbf{s}}
\newcommand\vecU{\mathbf{u}}

\newcommand\vecX{\mathbf{x}}
\newcommand\vecY{\mathbf{y}}

\newcommand\vecO{\boldsymbol{\omega}}

\newcommand\bbE{\mathbb{E}}
\newcommand\bbR{\mathbb{R}}
\newcommand\bbZ{\mathbb{Z}}

\newcommand{\iD}{d} 
\newcommand{\sS}{n} 

%% file: main.bbl
\begin{thebibliography}{22}
\providecommand{\natexlab}[1]{#1}
\providecommand{\url}[1]{\texttt{#1}}
\expandafter\ifx\csname urlstyle\endcsname\relax
  \providecommand{\doi}[1]{doi: #1}\else
  \providecommand{\doi}{doi: \begingroup \urlstyle{rm}\Url}\fi

\bibitem[{A. Zaytsev, E. Romanenkova, D. Ermilov}(2017)]{github}
{A. Zaytsev, E. Romanenkova, D. Ermilov}.
\newblock Interpolation errors for {G}aussian process regression.
\newblock
  \url{https://github.com/likzet/gp_interpolation_error/blob/master/code/gp_misspecification_expanded_experiments.ipynb},
  2017.

\bibitem[Bachoc(2013)]{bachoc2013cross}
F.~Bachoc.
\newblock Cross validation and maximum likelihood estimations of
  hyper-parameters of gaussian processes with model misspecification.
\newblock \emph{Computational Statistics \& Data Analysis}, 66:\penalty0
  55--69, 2013.

\bibitem[Bachoc(2018)]{bachoc2018asymptotic}
F.~Bachoc.
\newblock Asymptotic analysis of covariance parameter estimation for gaussian
  processes in the misspecified case.
\newblock \emph{Bernoulli}, 24\penalty0 (2):\penalty0 1531--1575, 2018.

\bibitem[Belyaev et~al.(2015)Belyaev, Burnaev, and
  Kapushev]{belyaev2015gaussian}
M.~Belyaev, E.~Burnaev, and Y.~Kapushev.
\newblock Gaussian process regression for structured data sets.
\newblock In \emph{International Symposium on Statistical Learning and Data
  Sciences}, pages 106--115. Springer, 2015.

\bibitem[Burnaev et~al.(2016)Burnaev, Panov, and
  Zaytsev]{burnaev2016regression}
E.V. Burnaev, M.E. Panov, and A.A. Zaytsev.
\newblock Regression on the basis of nonstationary {G}aussian processes with
  {B}ayesian regularization.
\newblock \emph{Journal of communications technology and electronics},
  61\penalty0 (6):\penalty0 661--671, 2016.

\bibitem[Castillo et~al.(2008)]{castillo2008lower}
I.~Castillo et~al.
\newblock Lower bounds for posterior rates with {G}aussian process priors.
\newblock \emph{Electronic Journal of Statistics}, 2:\penalty0 1281--1299,
  2008.

\bibitem[Cressie(2015)]{cressie2015statistics}
N.~Cressie.
\newblock \emph{Statistics for spatial data}.
\newblock John Wiley \& Sons, 2015.

\bibitem[Golubev and Krymova(2013)]{golubev13interpolation}
G.K. Golubev and E.A. Krymova.
\newblock On interpolation of smooth processes and functions.
\newblock \emph{Problems of Information Transmission}, 49\penalty0
  (2):\penalty0 127--148, 2013.

\bibitem[{GPy}(2012)]{gpy2014}
{GPy}.
\newblock {GPy}: {A {G}aussian process framework in python}.
\newblock \url{http://github.com/SheffieldML/GPy}, 2012.

\bibitem[Kolmogorov(1941)]{kolmogorov1992interpolation}
A.N. Kolmogorov.
\newblock Interpolation and extrapolation of stationary random sequences.
\newblock \emph{Izv. Akad. Nauk SSSR, Ser. Mat.}, 5\penalty0 (1):\penalty0
  3--14, 1941.

\bibitem[Le~Gratiet and Garnier(2015)]{le2015asymptotic}
L.~Le~Gratiet and J.~Garnier.
\newblock Asymptotic analysis of the learning curve for gaussian process
  regression.
\newblock \emph{Machine Learning}, 98\penalty0 (3):\penalty0 407--433, 2015.

\bibitem[Minasny and McBratney(2005)]{minasny2005matern}
B.~Minasny and A.~McBratney.
\newblock The {M}at{\'e}rn function as a general model for soil variograms.
\newblock \emph{Geoderma}, 128\penalty0 (3):\penalty0 192--207, 2005.

\bibitem[Panov(2016)]{panov2016nonasymptotic}
M.E. Panov.
\newblock Nonasymptotic approach to {B}ayesian semiparametric inference.
\newblock In \emph{Doklady Mathematics}, volume~93, pages 155--158. Springer,
  2016.

\bibitem[Rasmussen and Williams(2006)]{rasmussen2006}
C.~E. Rasmussen and C.~K.~I. Williams.
\newblock \emph{{G}aussian processes for machine learning}.
\newblock The MIT Press, 2006.

\bibitem[Stein(2012)]{stein2012interpolation}
M.~Stein.
\newblock \emph{Interpolation of spatial data: some theory for kriging}.
\newblock Springer Science \& Business Media, 2012.

\bibitem[Suzuki(2012)]{suzuki2012pac}
T.~Suzuki.
\newblock {PAC}-{B}ayesian bound for {G}aussian process regression and multiple
  kernel additive model.
\newblock In \emph{COLT}, pages 8--1, 2012.

\bibitem[Vaart and Zanten(2011)]{vaart2011information}
A.~van~der Vaart and H.~van Zanten.
\newblock Information rates of nonparametric gaussian process methods.
\newblock \emph{Journal of Machine Learning Research}, 12\penalty0
  (Jun):\penalty0 2095--2119, 2011.

\bibitem[van~der Vaart and van Zanten(2008)]{van2008rates}
A.~van~der Vaart and J.~van Zanten.
\newblock Rates of contraction of posterior distributions based on {G}aussian
  process priors.
\newblock \emph{The Annals of Statistics}, pages 1435--1463, 2008.

\bibitem[Wiener(1949)]{wiener1949extrapolation}
N.~Wiener.
\newblock \emph{Extrapolation, interpolation, and smoothing of stationary time
  series}, volume~2.
\newblock MIT press Cambridge, MA, 1949.

\bibitem[Wilcoxon(1945)]{wilcoxon1945individual}
F.~Wilcoxon.
\newblock Individual comparisons by ranking methods.
\newblock \emph{Biometrics bulletin}, 1\penalty0 (6):\penalty0 80--83, 1945.

\bibitem[Zaytsev and Burnaev(2017)]{zaytsev2017minimax}
A.~Zaytsev and E.~Burnaev.
\newblock Minimax approach to variable fidelity data interpolation.
\newblock In \emph{Artificial Intelligence and Statistics}, pages 652--661,
  2017.

\bibitem[Zaytsev et~al.(2014)Zaytsev, Burnaev, and
  Spokoiny]{zaytsev2014properties}
A.~Zaytsev, E.~Burnaev, and V.~Spokoiny.
\newblock Properties of the {B}ayesian parameter estimation of a regression
  based on {G}aussian processes.
\newblock \emph{Journal of Mathematical Sciences}, 203\penalty0 (6):\penalty0
  789--798, 2014.

\end{thebibliography}
